\documentclass[preprint,12pt]{elsarticle}
\usepackage{mathrsfs}
\usepackage{amsmath}
\usepackage{amsfonts}
\usepackage{amsthm}
\usepackage{amssymb}
\begin{document}

\begin{frontmatter}



\title{{\bf Global weak solutions to the two-dimensional Navier-Stokes equations of
compressible heat-conducting flows with symmetric data and forces}}


\author{Fei Jiang\corref{cor1}}
\ead{jiangfei0591@163.com}

\author{Song Jiang}
 \author{Junpin Yin}
 \cortext[cor1]{Corresponding author: Tel +86+15001201710.}
\address{Institute of Applied Physics and Computational Mathematics,
P.O. Box 8009, Beijing 100088, China.}

\begin{abstract}
We prove the global existence of weak solutions to the Navier-Stokes
equations of compressible heat-conducting fluids in two spatial
dimensions with initial data and external forces which are large and
spherically symmetric. The solutions will be obtained
as the limit of the approximate solutions in an annular domain.
We first derive a number of regularity results on the
approximate physical quantities in the ``fluid region'', as well as
the new uniform integrability of the velocity and temperature in the entire
space-time domain by exploiting the theory of the Orlicz spaces.
By virtue of these a priori estimates we then argue in a
manner similar to that in [Arch. Rational Mech. Anal. {\bf 173} (2004),
297-343] to pass to the limit and show that the limiting functions are
indeed a weak solution which satisfies the mass and momentum
equations in the entire space-time domain in the sense of
distributions, and the energy equation in any compact
subset of the ``fluid region''.
\end{abstract}


\end{frontmatter}


\newtheorem{thm}{Theorem}[section]
\newtheorem{lem}{Lemma}[section]
\newtheorem{pro}{Proposition}[section]
\newtheorem{cor}{Corollary}[section]
\newproof{pf}{Proof}
\newdefinition{rem}{Remark}[section]
\newtheorem{definition}{Definition}[section]

\section{Introduction}\numberwithin{equation}{section}
The two-dimensional Navier-Stokes equations of compressible heat-conducting fluids
express the conservation of mass, and the balance of momentum and
energy, which can be written as follows in Eulerian coordinates.
\begin{eqnarray} \label{0101}  &&
\varrho_t+\mathrm{div}(\varrho\mathbf{u})=0,  \\
&&  \label{0102}
(\varrho\mathbf{u})_t+\mathrm{div}(\varrho\mathbf{u}\otimes\mathbf{u})+\nabla P(\varrho,\theta)
=\mu \Delta\mathbf{u}+(\lambda+\mu)\nabla\mathrm{div}\mathbf{u} +\varrho\mathbf{f},  \\
&& \label{0103}
(\varrho E)_t+\mathrm{div}(\varrho E+P(\varrho,\theta)\mathbf{u}) =
\Delta\Big(\kappa\theta+\frac{1}{2}\mu|\mathbf{u}|^2\Big)+\mu
\mathrm{div}\big[(\nabla\mathbf{u})\mathbf{u}\big]  \nonumber \\
&& \qquad \qquad \quad  \qquad\qquad \qquad \qquad\;
+\lambda\mathrm{div}\big[(\mathrm{div}\mathbf{u})\mathbf{u}\big]
+\varrho\mathbf{u}\cdot\mathbf{f}.
\end{eqnarray}
Here $\varrho$, $\mathbf{u}=(u_1,u_2)\in\mathbb{R}^2$, and $\theta$,
$$E=\frac{|u|^2}{2} +\theta\qquad \mbox{and}\qquad P=K\varrho\theta$$
are the density, velocity,
temperature, total energy density and pressure of an ideal gas (with
unit specific heat), respectively; $\mu$ and $\lambda$ are the constant viscosity coefficients
satisfying $\mu >0$ and $\mu+\lambda\geq 0$, $\kappa >0$ is the heat-conduction coefficient;
$\mathbf{f}=(f_1,f_2)$ is the external force, and $\nabla \mathbf{u}$
denotes the gradient of the velocity vector with respect to the spatial variable
$\mathbf{x}\in {\mathbb R}^2$.

We consider an initial boundary value problem for the system (\ref{0101})--(\ref{0103}) in
a ball $\Omega :=\{\mathbf{x}\in\mathbb{R}^2;\; |\mathbf{x}|<R\}$ with boundary conditions
\begin{equation}\label{0104}
\mathbf{u}=\mathbf{0},\quad \frac{\partial\theta}{\partial\mathbf{n}}=0\;\;\mbox{ on }\partial\Omega
\end{equation}
and initial conditions
\begin{equation}\label{0105}
(\varrho,\mathbf{u},\theta)|_{t=0}=(\varrho,\mathbf{u}_0,\theta_0),
\end{equation}
where $\mathbf{n}$ is the outer normal vector to $\partial\Omega$.

 In the spherically symmetric case, namely,

\begin{equation}\label{0106}
\varrho(t,\mathbf{x})=\varrho(t,r),\;\;
\mathbf{u}(t,\mathbf{x})=u(t,r)\frac{\mathbf{x}}{r},\;\; \theta(t,\mathbf{x})=\theta (t,r),\;\;
\mathbf{f}(t,\mathbf{x})=f(t,r)\frac{\mathbf{x}}{r},
\end{equation}
where $r=r(\mathbf{x}):=|\mathbf{x}|$,
the system (\ref{0101})-(\ref{0103}) takes the form:
\begin{eqnarray}
&& \varrho_t+(\varrho u)_r+\frac{\varrho u}{r}=0, \label{0107} \\
&& (\varrho u)_t+(\varrho u^2)_r+\frac{\varrho u^2}{r}+P_r(\varrho,\theta)
-\nu \left(u_r+\frac{u}{r}\right)_r=\varrho f,  \label{0108} \\
&& (\varrho\theta)_t+(\varrho u \theta)_r+\frac{\varrho u
\theta}{r}-\kappa\left(\theta_{rr}+\frac{\theta_r}{r}\right)+P(\varrho,\theta)
\left(u_r+\frac{u}{r}\right)=\mathcal{Q},  \label{0109} \end{eqnarray}
where
\begin{equation}\label{0110}
\nu:=\lambda+2\mu ,\quad\;\;\mathcal{Q}=\nu\left(u_r+\frac{u}{r}\right)^2 -\frac{2\mu}{r}
\partial_ru^2\geq 0.\end{equation}
The boundary and initial conditions become
\begin{eqnarray}\label{01101} &&
u=\theta_r =0\qquad\mbox{at }r=R , \label{0109a} \\[0.5em]
&& (\varrho,u,\theta)|_{t=0}=(\varrho,u_0,\theta_0), \qquad 0<r < R.  \label{0109b}
\end{eqnarray}

The main purpose of this paper is to prove the global existence of
weak solutions to the problem (\ref{0101})--(\ref{0105}) when the
initial data and external forces are large and spherically
symmetric. Our work is motivated by the paper of Hoff and Jenssen
\cite{HDJHK} where they studied the spherically and cylindrically
symmetric nonbarotropic flows with large data and forces, and
established the global existence of weak solutions to the
compressible nonbarotropic Navier-Stokes equations in the ``fluid
region''. In the entire space-time domain, however, the momentum
equation in \cite[Thorem 1.1]{HDJHK} only holds weakly with a
nonstandard interpretation of the viscosity terms as distributions.
A natural question is to ask whether the momentum equation holds in
the standard sense of distributions. A positive answer was given
recently by Zhang, Jiang and Xie \cite{ZJWJSXF} for a screw pinch
model arisen from plasma physics when the heat conductivity $\kappa$
satisfies certain growth conditions. In the present paper, based on
some new uniform global estimates of $u$ and $\theta$ (see Lemma
\ref{lem:0204}) which are established by applying the theory of the
Orlicz spaces, we can give a positive result for the two-dimensional
Navier-Stokes equations of compressible heat-conducting fluids
(\ref{0101})--(\ref{0103}), improving therefore the result of
\cite{HDJHK}.

We now give the precise statement of our assumptions and results. The
external force $\mathbf{f}$ is assumed to satisfy
\begin{equation}\label{0111}
\mathbf{f}\in L^1((0,T), L^\infty(\Omega))\cap
L^\infty((0,T),L^2(\Omega))\end{equation} for each $T>0$. The
initial data ($\varrho_0,\mathbf{u}_0,\theta_0$) are assumed to
satisfy
\begin{equation}\label{0112}
C_0^{-1}\leq \varrho_0\leq C_0,\ C_0^{-1}\leq \theta_0\ \mbox{a.e.
in }\Omega,\end{equation}
\begin{equation}\label{0113}
\int_\Omega\varrho_0S(\varrho_0,\mathbf{u}_0,\theta_0)\mathrm{d}\mathbf{x}\leq
C_0\qquad \mbox{for some positive constant }C_0,  \end{equation}
where $S$ is the entropy density in the form of
\begin{equation}\label{0114}
S(\varrho,\mathbf{u},\theta)=K\Psi(\varrho^{-1})+\Psi(\theta)+\frac{1}{2}|\mathbf{u}|^2
\end{equation}
with $\Psi(s) = s -\mathrm{log}s -1$. We point out here that there are no
smallness or regularity conditions imposed on $\mathbf{f}$ and
$(\varrho_0,\mathbf{u}_0,\theta_0)$.

Under the conditions (\ref{0111})--(\ref{0114}), we shall prove the
following existence theorem on spherically symmetric solutions
to the problem (\ref{0101})--(\ref{0105}).
\begin{thm}\label{thm:0101}
Assume that the initial data ($\varrho_0,\mathbf{u}_0,\theta_0$) and
the external force $\mathbf{f}$ are spherically symmetric and
 satisfy the conditions (\ref{0111})--(\ref{0113}). Then, the initial boundary
problem (\ref{0101})--(\ref{0105}) has a global weak solution
($\varrho,\mathbf{u},\theta$) in the form of (\ref{0106}) satisfying the following:
\begin{enumerate} [\quad \ (a)]
  \item  The support of $\varrho$ is
bounded on the left by  a H$\ddot{o}$lder curve $\underline{r}(t)\in
C^{0,1/4}_{\mathrm{loc}}:[0,\infty)\rightarrow [0,\infty)$.
Moreover, if $\mathcal{F}$ is the ``fluid region", defined by
\begin{equation*}\label{0115}
\mathcal{F}:=\{(t,\mathbf{x})~|~t\geq 0\ and\
\underline{r}(t)<r(t)\leq R\},\end{equation*}  then $\mathcal{F}\cap
\{t>0\}\cap \{r<R\}$ is open set.
  \item The density $\varrho\in
  L^\infty_{\mathrm{loc}}(\mathcal{F})$, $\mathbf{u}$ and $\theta$
  are locally H\"{o}lder continuous in $\mathcal{F}\cap\{t>0\}$,
  and the Navier-Stokes equations
  (\ref{0101})--(\ref{0103}) hold in $\mathcal{D}'(\mathcal{F}\cap\{t>0\}\cap
  \{r<R\})$.
\item The density $\varrho \in
  C([0,\infty),W^{1,\infty}(\Omega)^*)$. Also, $\varrho(t,\cdot)\equiv
  0$ in $\Omega\backslash\bar{\mathcal{F}}$, and if $\varrho
  \mathbf{u}$ is taken to be zero in
  $\Omega\backslash\bar{\mathcal{F}}$, then the weak form of the
  mass equation (\ref{0101}) holds for test functions $\psi\in C^1([t_1,t_2]\times
  \bar{\Omega})$:
\begin{equation}\label{0116}
\int_\Omega \varrho \psi
\mathrm{d}\mathbf{x}\bigg|_{t_1}^{t_2}=\int_{t_1}^{t_2}\int_\Omega
(\varrho \psi_t+\varrho\mathbf{u}\cdot\nabla\psi)\mathrm{d}\mathbf{x}\mathrm{d}t.
\end{equation}
  \item The velocity $\mathbf{u}\in L^{4/3}_{\mathrm{loc}}( [0,\infty),W^{1,4/3}(\Omega))$.
  For $t_1<t_2$, if $\psi\in C^1([t_1,t_2]\times
  \bar{\Omega})$ vanishes on $\partial \Omega$, then, for $i=1,2$,
\begin{eqnarray}
&& \int_\Omega\varrho{u}_i\psi\mathrm{d}\mathbf{x}\bigg|_{t_1}^{t_2}-\int_{t_1}^{t_2}\int_\Omega
\left(\varrho u_i\psi_t +\varrho{u}_i\mathbf{u}\cdot \nabla \psi +
P(\varrho,\theta)\psi_{x_i}\right)\mathrm{d}\mathbf{x}\mathrm{d}t\nonumber \\
= && \int_{t_1}^{t_2}\int_\Omega\varrho
f_i\mathrm{d}\mathbf{x}\mathrm{d}t -\int_{t_1}^{t_2}\int_\Omega
[(\lambda+\mu){\mathbf{u}_{x_i}}\cdot
\nabla\psi+\mu\nabla{u}_i\nabla\psi]\mathrm{d}\mathbf{x}\mathrm{d}t.
\label{0117}
\end{eqnarray}
  \item The gradient
  $\nabla$$\theta$$\in$$L^1_{\mathrm{loc}}(\mathcal{F})$, and the weak form of the energy
  equation (\ref{0103}) holds for test functions $\psi\in C^1([t_1,t_2]\times\bar\Omega)$
  for which there is an $\eta >0$ such that {\rm supp}$\psi(t,\cdot)\subset
   \{\mathbf{x}~|~\underline{r}(t)+\eta\leq
  r(\mathbf{x})\}$ for each $t\in [t_1,t_2]$:
\begin{eqnarray}
&& \int_\Omega \varrho
E\psi\mathrm{d}\mathbf{x}\bigg|_{t_1}^{t_2}-\int_{t_1}^{t_2}\int_\Omega(\varrho
E\psi_t+(\varrho E+P(\varrho,\theta))\mathbf{u}\cdot \nabla
\psi)\mathrm{d}\mathbf{x}\mathrm{d}t = \int_{t_1}^{t_2}\int_\Omega\varrho \mathbf{f}\cdot
\mathbf{u}\psi\mathrm{d}\mathbf{x}\mathrm{d}t \qquad \nonumber \\
&& \quad - \int_{t_1}^{t_2}\int_\Omega\left(\kappa\nabla \theta
+\frac{1}{2}\mu\nabla |\mathbf{u}|^2 +\mu (\mathrm{div}\mathbf{u})\mathbf{u}
+\lambda (\nabla\mathbf{u})\mathbf{u}\right)\cdot\nabla \psi\mathrm{d}\mathbf{x}\mathrm{d}t.
\label{0118}
\end{eqnarray}
  \item The total energy energy, minus the mechanical work done by
  the external force, is weakly nonincreasing in time. That is, if
\begin{equation*}\label{0119}
\mathcal{E}(t):=\int_\Omega\varrho(t,\mathbf{x})\left[\theta(t,\mathbf{x})+\frac{1}{2}|\mathbf{u}|^2
\right]\mathrm{d}\mathbf{x}\mathrm{d}t,
\end{equation*}
 then
\begin{equation*}\label{0120}
\mathcal{E}(t)=\mathcal{E}(0)+\int_0^t\int_\Omega \varrho
\mathbf{u}\cdot
\mathbf{f}\mathrm{d}\mathbf{x}\mathrm{d}t-\lim_{b\rightarrow 0}
\lim_{j\rightarrow \infty}\int_{\varepsilon_j\leq r(\mathbf{x})\leq b}
(\varrho^j E^j)(t,\mathbf{x})\mathrm{d}\mathbf{x}
\end{equation*}
as a function of $t$ in $\mathcal{D}'(0,\infty)$, where $E^j=|\mathbf{u}^j|^2/2 + \theta^j$.
\end{enumerate}
\end{thm}

\begin{rem} In \cite{HDJHK}, Hoff and Jenssen proved
that (\ref{0102}) holds in $\mathcal{D}'(\mathcal{F}\cap\{ t>0\}\cap
\{r<R\})$. Here, by exploiting the theory of the Orlicz spaces we
are able to derive some new uniform global integrability of the
approximate solutions (cf. Lemma \ref{lem:0204}) to show that
(\ref{0102}) holds in the entire space-time domain in the weak sense
(i.e., (d) of Theorem 1.1). Furthermore, this method can be applied
to the screw pinch model with positive  constant heat-conduction
coefficient in \cite{ZJWJSXF} and the
 cylindrically symmetric rotating model of (\ref{0101})--(\ref{0103})
 (that is, in the symmetric equations (6)--(10) in \cite{HDJHK}, we take
  $u=v$, $w=0$, $f_1=f_2$ and $f_3=0$) to obtain similar results.\end{rem}

Combining the global a priori estimates derived in Subsection 2.3,
we shall prove Theorem \ref{0101} in Section 4 by the convergence
argument similar to that in \cite{HDJHK}. For this purpose, we
consider the approximate solutions
($\varrho^j,\mathbf{u}^j,\theta^j$) of the problem
(\ref{0101})--(\ref{0105}) in the annular regions
$\Omega^j:=\{\mathbf{x}~|~\varepsilon_j<r(\mathbf{x})<R\}$, where
$\varepsilon_j$ is a sequence of positive inner radii tending to 0.
Since the $1/r$ singularity in the equations
(\ref{0107})--(\ref{0109}) plays no role at the stage when
$\varepsilon_j$ is fixed and positive, the global existence of
approximate solutions ($\varrho^j,\mathbf{u}^j,\theta^j$) for
(\ref{0101})--(\ref{0105}) can thus be shown in a manner similar to
that in \cite{HDJHK,ZJWJSXF}. However, to pass to the limit as
$j\rightarrow \infty$ and to show the global existence of weak
solutions to the original problem (\ref{0101})--(\ref{0105}), we need
some $\varepsilon_j$-independent a priori estimates. This will be
done in Sections 2 and 3. We first prove the global estimates in
Section 2, where we derive the standard energy-entropy estimates in
Subsection 2.1, and apply these estimates to establish a new uniform
integrability of the approximate solutions in the entire spacetime
domain by exploiting the theory of Orlicz spaces in Subsections 2.2
and 2.3, which is crucial in the proof of (d) of Theorem 1.1. Then,
in Section 3, we list the well-known the pointwise bounds for
$\varrho_j$ and $\theta_j$ as consequences of the energy and entropy
estimates. These pointwise bounds are independent of
$\varepsilon_j$, but only away from the origin of Lagrangian space.
More precisely, as in \cite{HDJHK}, for any given $h>0$ we define
the particle position $r_h^j(t)$ by
\begin{equation*}\label{0121}
h=\int_{\varepsilon_j}^{r_h^j(t)}\varrho^j(t,r)r\mathrm{d}r
\end{equation*}
from which and the standard energy-entropy estimates it follows that
there exists a positive constant $C(h)$, depending only on $h>0$,
such that $r_h^j(t)\geq C(h)> 0$. With this observation, we can obtain that
for any fixed $h>0$ and $T>0$, there is a positive
constant $C(T,h)$, depending only on $h$, $T$ and the initial data, such that
\begin{equation*}\label{0123}
C(T,h)^{-1}\leq \varrho^j(t,\mathbf{x})\leq C(T,h)\;\;\mbox{ for any
}(t,\mathbf{x})\in [0,T]\times[r_h^j(t),R].
\end{equation*}
 Applying these pointwise bounds, we can get a number of
higher-order energy estimates for the approximate solutions in
Subsection 3.2, which are also independent of $\varepsilon_j$ and
only away from the origin of Lagrangian space. These
$\varepsilon_j$-independent bounds enable us to define the ``fluid
region" $\mathcal{F}$ (see (a) of Theorem 1.1) and to obtain the
uniform H\"{o}lder continuity of the quantities on the compact
subsets of $\mathcal{F}\cap \{t>0\}$ (see (b) of Theorem 1.1).
Finally, all the assertions of (a)--(f) indicated in Theorem 1.1
will be proved in Section 4 by the convergence arguments adapted
from Hoff and Jenssen's paper \cite{HDJHK}.  We note that the final
step of this argument provides a sort of a posteriori validation
that the equations (\ref{0107})--(\ref{0109}) are indeed the correct
forms of the general system (\ref{0101})--(\ref{0103}) in the
symmetric case considered here.

As pointed out in \cite{HDJHK}, we still do not have sufficient
information to infer that $\underline{r}(t)\equiv 0$, nor do we know
whether solutions exist for which $\underline{r}\neq 0$. The
analysis simply shows that $\underline{r}(t)$ may be positive, and
that, if it is, a vacuum state of radius $\underline{r}(t)$ centered
at the origin. In any case, the total mass is conserved in the
spherical case, as is clear from (c) of Theorem 1.1, and the total
momentum is zero because of symmetry.

We show in (e) only that the energy equation holds on the support of
$\varrho$, rather than in the entire space-time domain
$(0,\infty)\times \Omega$. This is partly due to that we cannot
obtain higher global regularity of $\theta$ and $\mathbf{u}$. We may
regard the restriction in (5) that the test function be supported in
$\mathcal{F}$ as reasonable, since there is no fluid outside
$\mathcal{F}$, and the model is not really valid there.
Additionally, the failure of the analysis to detect whether or not
energy is lost ((f) of Theorem \ref{thm:0101}) calls into question
the adequacy of the mass, energy, and entropy bounds in Lemma 2.1,
which are the only known (global) a priori bounds in the
multidimensional case now.

We end this section by mentioning some related existence results for
large data in the multidimensional case. The global existence of
weak solutions was first shown by Lions \cite{LPLMTFM98} for
isentropic flows under the assumption that the specific heat ratio
$\gamma > 3n/(n+2)$ where $n=2,3$ denotes the spatial dimension.
Then, by using the curl-div lemma to delicately derive certain
compactness, and applying Lions' idea and a technique from
\cite{JZ01}, Feireisl, Novotn\'y and Petzeltov\'a \cite{Fe01,FNP01}
extended Lions' existence result to the case $\gamma > n/2$. For any
$1\leq \gamma \leq n/2$, a global weak solution still exists when
the initial data have certain symmetry (e.g., spherical, or
axisymmetric symmetry), see \cite{HDSI}, \cite{JZ03}--\cite{JZ01}.
For non-isentropic flows, the global existence for general data is
still not available. Recently, under certain growth conditions upon
the pressure, viscosity and heat-conductivity (i.e., radiative
gases), Feireisl, et al. obtained the global existence of the
so-called ``variational solutions'' in the sense that the energy
equation is replaced by an energy inequality, see \cite{FETD2} for
example. However, this result excludes the case of ideal gases
unfortunately. The global existence of a solution for large data in
the non-isentropic case needs further study.

\section{Global Estimates}
In this section we derive a priori global estimates for any smooth
(approximate) solution ($\varrho^\varepsilon$, $u^\varepsilon$,
$\theta^\varepsilon$) of (\ref{0107})--(\ref{0109}) together with
additional boundary conditions:   \begin{equation}
u^\varepsilon=\theta_r^\varepsilon=0\qquad\mbox{at
}r=\varepsilon\mbox{ and }R. \label{j1}  \end{equation}

 We assume that the initial data and force are
smooth and satisfy the bounds (\ref{0111})--(\ref{0113}) with
constants independent of $\varepsilon$ and
$$ \int_\varepsilon^R\varrho_0^\varepsilon\mathrm{d}\mathbf{x}\equiv
M_0:=\int_0^R\varrho_0\mathrm{d}\mathbf{x} .$$
We refer to Section 4 for a brief discussion on the existence of such approximate
solutions. As discussed in Section 1, we shall eventually take a
sequence of inner radii $\varepsilon_j\rightarrow 0$ to prove Theorem 1.1. Since
$\varepsilon>0$ is fixed for the time being, we suppress the
dependence on $j$.

\subsection{Energy and Entropy Estimates}
We start with the following lemma which states the standard energy and entropy estimates
for these approximate solutions.
\begin{lem}\label{lem:0201}
Let ($\varrho^\varepsilon$, $u^\varepsilon$, $\theta^\varepsilon$)
be a smooth solution of (\ref{0107})--(\ref{0109b}) defined on
$[0,T]\times [\varepsilon, R]$ with boundary conditions (\ref{j1}).
Then, there are constants $M_0$ and $C(T)$, such that
\begin{eqnarray}\label{0201}
&& \int_\varepsilon^R \varrho^\varepsilon(t,r)r\mathrm{d}r\equiv M_0, \\
&& \label{0202}
\int_\varepsilon^R (\varrho^\varepsilon E^\varepsilon)(t,r)r
\mathrm{d}r\leq C(T), \\
&& \label{0203}  \int_\varepsilon^R (\varrho^\varepsilon
S^\varepsilon)(t,r)r\mathrm{d}r+\int_0^t\int_\varepsilon^R
\left[\kappa\left(\frac{\theta^\varepsilon_r}{\theta^\varepsilon}\right)^2
+\frac{\mathcal{Q}}{\theta^\varepsilon} \right]r\mathrm{d}r\leq C(T)
\quad \mbox{for all }t\in [0,T],  \end{eqnarray} where
$E^\varepsilon=(u^\varepsilon)^2/{2}+\theta^\varepsilon$,
$S^\varepsilon$ is the entropy density defined in (\ref{0114}) and
$\mathcal{Q}$ is given in (\ref{0110}).
\end{lem}
\begin{pf}
The bounds (\ref{0201})--(\ref{0203}) are the standard energy estimates which
follow directly from the equations (\ref{0107})--(\ref{0109}),
the boundary conditions and the assumption (\ref{0111}) on the external force.
\hfill $\square$
\end{pf}
\subsection{Excursion to Theory of the Orlicz Spaces}
Before deriving the global estimates on the temperature and velocity, we
recall some well-known results concerning the Orlicz
spaces (see, for example, \cite{ARAJJFF,KAJOFS} for details), which are often
used to investigate the 2D compressible Navier-Stokes equations (see
\cite{JSZPRI,JFTZOM32,EROM} for example).
\begin{definition}[Young's function] We say that $\Phi$ is a
Young's (or $N$-) function if
\begin{equation*}\label{0204}
\Phi(t)=\int_0^t\phi(s)\mathrm{d}s,\quad t\geq 0,
\end{equation*}
where the real-valued function $\phi$ defined on $[0,\infty)$ has
the following properties
\begin{equation*}\label{0205}\begin{array}{ll}
\phi(0)=0,\;\;\phi(s)>0,\ s>0,\;\;\lim_{s\rightarrow \infty}
\phi(s)=\infty, \\[0.5em] \phi \mbox{ is right continuous and
nodecreasing on }[0,\infty).\end{array}\end{equation*} We define
\begin{equation*}\label{0206}\psi(t)=\sup_{\{\phi(s)\leq t\}}s,\;\; t\geq 0,\quad
\Psi(t)=\int_0^t\psi(s)\mathrm{d}s.  \end{equation*} Then $\Psi$ is
a Young's function as well. We call $\Psi$ the complementary Young's
function to $\Phi$. If $\Phi$ is complimentary to $\Psi$, then
$\Psi$ is complimentary to $\Phi$.
\end{definition}
\begin{definition}[Orlicz spaces] Let $\Omega$ be a domain in
$\mathbb{R}^n$ and let $\Phi$ be a Young function. The Orlicz class
$K_\Phi(\Omega)$ is the set of all (equivalent classes modulo
equality a.e. in $\Omega$ of) measure functions $u$ defined on
$\Omega$ that satisfy $\int_\Omega
\Phi(|u(\mathbf{x})|)\mathrm{d}\mathbf{x}<\infty$. The Orlicz space
$L_\Phi(\Omega)$ is the linear hull of the Orlicz class
$K_\Phi(\Omega)$, that is, the smallest vector space that contains
$K_\Phi(\Omega)$. The functional
\begin{equation*}\label{0207}\|u\|_{\Phi(\Omega)}=\inf\left\{k>0~\bigg|
\int_\Omega\Phi\left(\frac{|u(\mathbf{x})|}{k}\right)
\mathrm{d}\mathbf{x}\leq1\right\}<\infty\end{equation*} is a norm on
$L_\Phi(\Omega)$. It is called the Luxembourg norm. Thus,
$L_\Phi(\Omega)$ is a Banach space with respect to the Luxembourg
norm.
  \end{definition}
\begin{definition}[Cone condition] Let $\mathbf{y}$ be
a nonzero vector in $\mathbb{R}^n$. Let $\angle
(\mathbf{x},\mathbf{y})$ be the angle between the position vector
$\mathbf{x}$ and $\mathbf{y}$. For given such $\mathbf{y}$, $h>0$,
and $k$ satisfying $0<k\leq \pi$, the set
\begin{equation*}
\Lambda=\{\mathbf{x}\in \mathbb{R}^n~|~\mathbf{x}=0\mbox{ or
}0<|\mathbf{x}|\leq h,\ \angle (\mathbf{x},\mathbf{y})\leq k/2\}
\end{equation*}
is called a finite cone of height $h$, axis direction $\mathbf{y}$
and aperture angle $k$ with vertex at the origin.

$\Omega\subset \mathbb{R}^n$ satisfies the cone condition if there exists a finite
cone $\Lambda$, such that each $\mathbf{x}\in \Omega$ is the vertex of a finite
cone $\Lambda_{\mathbf{x}}$ contained in $\Omega$ and congruent to $\Lambda$.
\end{definition}

Now, we define
$$M=M(s):=(1+s)\mathrm{ln}(1+s)-s,\;\; N=N(s):=e^s-s-1,\;\; H=H(s):=e^{s^2}-1.$$
Next, we list some basic facts on the Orlicz spaces
$L_M(\Omega)$, $L_N(\Omega)$ and $L_H(\Omega)$.
\begin{enumerate}[\quad \ (a)]
\item $M$ and $N$ are the complementary Young's functions (see \cite[8.3]{ARAJJFF}).
 \item Let $u(\mathbf{x})\in L_{M}(\Omega)$ and $v(\mathbf{x})\in L_{N}(\Omega)$.
By virtue of the generalized H\"{o}lder's inequality
(see \cite[8.11]{ARAJJFF}), we have  $uv\in L^1(\Omega)$ and
\begin{equation}\label{0208}
\left|\int_\Omega uv\mathrm{d}\mathbf{x}\right|\leq
2\|u\|_{M(\Omega)}\|v\|_{N(\Omega)}.
\end{equation}
\item  Let  $\Omega$  be bounded
and satisfy the cone condition in $\mathbb{R}^2$. By virtue of
\cite[Theorem 8.12, 8.25 and 8.27]{ARAJJFF}, we have for any $p\geq
1$ that
\begin{equation}\label{0209}W^{1,2}(\Omega)\hookrightarrow L_{H}(\Omega)\hookrightarrow
L_{N}(\Omega)\hookrightarrow L^p(\Omega)\;\mbox{ and }\;
W^{1,2}(\Omega)\hookrightarrow\hookrightarrow
L_{N}(\Omega).\end{equation}
  \item Denote by $E_N(\Omega)$ the closure of the set of all bounded measurable
functions on $\Omega$ with respect to the Luxembourg norm
$\|\cdot\|_{N(\Omega)}$. Then, the Orlicz space $L_M(\Omega)$ is $E_N$-weakly compact,
i.e., for any sequence $\{v_n\}\in L_{M}(\Omega)$ uniformly bounded, there is a
  subsequence of $\{v_n\}$, still denoted by $\{v_n\}$, and a $v\in L_M(\Omega)$,
  such that
  $$\int_\Omega v_n\varphi\mathrm{d}\mathbf{x}\rightarrow \int_\Omega
  v\varphi\mathrm{d}\mathbf{x}\quad\;\mbox{for any }\varphi\in E_N(\Omega) $$
(see \cite[3. Appendix]{JSZPRI}).
\end{enumerate}

\subsection{Global Estimates on the Temperature and Velocity}
Now, we are in a position derive the global estimates on temperature and velocity.
First, we define
\begin{equation*}\label{0210}
\tilde{\varrho^\varepsilon}(\cdot,|\mathbf{x}|)=\left\{
                              \begin{array}{ll}
   \varrho^\varepsilon(\cdot,|\mathbf{x}|), &  \varepsilon<|\mathbf{x}|<R, \\
                                0, &
 0\leq |\mathbf{x} |\leq\varepsilon,
                              \end{array}
       \right.\ \tilde{\theta^\varepsilon}(\cdot,|\mathbf{x}|)=\left\{
        \begin{array}{ll}
         \theta^\varepsilon(\cdot,|\mathbf{x}|), &  \varepsilon<|\mathbf{x}|<R, \\
                                \theta^\varepsilon(\cdot,\varepsilon), &
 0\leq |\mathbf{x} |\leq\varepsilon ,
                              \end{array}
                            \right.
\end{equation*}
and make use of (\ref{0201})--(\ref{0203}) and the definition of
Luxemburg norm $\|\cdot\|_{L_M(\Omega)}$ to deduce that there exists a constant $C_1(T)$, such that
\begin{eqnarray} &&
\|\tilde{\varrho^\varepsilon}(t,|\mathbf{x}|)\|_{L^1(\Omega)}=M_0,\
\|\tilde{\varrho^\varepsilon}(t,|\mathbf{x}|)\|_{L_M(\Omega)}\leq
{C}_1(T),\qquad\forall\, t\in (0,T),    \label{0211} \\[1mm]
&& \int_{\Omega}
\tilde{\varrho^\varepsilon}(t,|\mathbf{x}|)\mathrm{ln}
\left(1+{\tilde{\theta^\varepsilon}(t,|\mathbf{x}|)}\right)\mathrm{d}\mathbf{x}\leq
{C}_1(T),\qquad\forall\, t\in (0,T),   \label{0212} \end{eqnarray}
and
\begin{equation}\label{0213}
\int_0^T\int_{\Omega}\left|\nabla
\mathrm{ln}\left(1+{\tilde{\theta^\varepsilon}(t,\mathbf{x})}\right)
\right|^2\mathrm{d}\mathbf{x}\mathrm{d}t\leq C_1(T).
\end{equation}

Notice that $\theta^\varepsilon(t,r)$ is a smooth function in
$(\varepsilon,R)$, we can easily verify that
\begin{equation}\label{0214}\mathrm{ln}(1+{\tilde{\theta^\varepsilon}(t,|\mathbf{x}|)
)}\in L^2(0,T; W^{1,2}(\Omega)).\end{equation} Furthermore, by
(\ref{0211}) and $L_{M}\hookrightarrow L^1(\Omega)$, we infer
that there exists a constant $C_2(\Omega)$, such that
\begin{equation}\label{0215}
C_2(\Omega)\leq
\|\tilde{\varrho^\varepsilon}(t,|\mathbf{x}|)\|_{L_M(\Omega)}\quad \mbox{for any }t\in (0,T).
\end{equation}

At this stage we shall need two auxiliary results: 1) The first one
is a revised version generalized Korn-Poincar\'{e} inequality (see
\cite[Theorem 10.17]{FENA}) in the case of the Orilcz spaces. The
idea of the proof is essentially the same as that used in the proof
of \cite[Lemma 3.2]{FETD2} under trivial modification. 2) The other
one is the revised version Sobolev embedding in two dimensions which
will be used to derive bounds of the temperature. This is an idea
due to Lions who ever used similar embedding to derive the global
integrability of the temperature in the proof of the existence of
weak solutions to the stationary problems for the full compressible
Navier-Stokes equations in a bounded domain
$\Omega\subset\mathbb{R}^2$ \cite[(6.204) in Section 6.11]{LPLMTFM98}. These two
 auxiliary results are formulated in the following two lemmas.
\begin{lem}[Generalized Korn-Poincar\'e inequality]\label{lem:0102}
Let $\Omega\subset\mathbb{R}^2$ be a bounded domain satisfying the cone condition. Assume
that $v\in W^{1,2}(\Omega)$, and $\varrho\geq 0$ satisfies
\begin{equation}\label{0216}
0< C_1\leq \|\varrho\|_{L(\Omega)},\ \|\varrho\|_{L_M(\Omega)}\leq C_2.
\end{equation}
Then there is a constant $C_3$ depending solely on $C_1$ and $C_2$, such that
\begin{equation*}\label{0217}
\left\|v\right\|_{L^2(\Omega)}\leq C_3(C_1, C_2)\left[ \left\|\nabla
v\right\|_{L^2(\Omega)}+\int_\Omega \varrho |v|\mathrm{d}\mathbf{x}
\right].
\end{equation*}
\end{lem}
\begin{pf}
We prove the lemma by contradiction. Suppose that
the conclusion of Lemma 2.2 be false, then there would be a sequence
$\{\varrho_n\}_{n=1}^\infty$ of non-negative functions satisfying
(\ref{0216}) and a sequence $\{v_n\}_{n=1}^\infty \subset W^{1,2}(\Omega)$, such that
\begin{equation}\label{0218}
\|v_n\|_{L^2(\Omega))}\geq C_n\left( \|\nabla
v_n\|_{L^2(\Omega)}+\int_\Omega \varrho_n |v_n|\mathrm{d}\mathbf{x}
\right),\quad C_n\rightarrow +\infty.
\end{equation}
Setting $w_n=v_n\|v_n\|_{L^2(\Omega)}^{-1}$, making use of
(\ref{0209}) and (\ref{0218}), we find that
\begin{equation}\label{0219}
w_n\rightarrow w=\frac{1}{\sqrt{\Omega}}\quad \mbox{ strongly in }L_N(\Omega).
\end{equation}

  In view of the hypothesis (\ref{0216}), we see that
\begin{equation}\label{0220}
\int_\Omega\varrho_n \varphi\mathrm{d}\mathbf{x}\to \int_\Omega\varrho
\varphi\mathrm{d}\mathbf{x}\quad \mbox{ for any }\varphi\in E_N(\Omega).
\end{equation}
Thus, by virtue of (\ref{0219}), (\ref{0220}) and (\ref{0208}), one has
\begin{equation}\label{0221}
\lim_{n\rightarrow  \infty}\int_\Omega\left(\varrho_n w_n-\varrho
w\right)\mathrm{d}\mathbf{x}=\lim_{n\rightarrow
\infty}\int_\Omega\varrho_n (w_n -w)\mathrm{d}\mathbf{x}+\lim_{n\rightarrow
\infty}\int_\Omega(\varrho_n -\varrho )w\mathrm{d}\mathbf{x}=0.
\end{equation}
The identity (\ref{0221}), together with (\ref{0216}) and (\ref{0219}), yields
\begin{equation}\label{0222}
\lim_{n\rightarrow  \infty}\int_\Omega\varrho_n w_n
\mathrm{d}\mathbf{x}=\int_\Omega\varrho w\mathrm{d}\mathbf{x}>0.
\end{equation}

 On the other hand, (\ref{0218}) implies
\begin{equation}\label{0223}
\lim_{n\rightarrow \infty}\int_\Omega\varrho_n w_n\mathrm{d}\mathbf{x}=0,
\end{equation}
which contradicts with (\ref{0222}). Therefore, the conclusion of Lemma 2.2 remains true.
\hfill $\square$
\end{pf}
\begin{lem}[Sobolev embedding]\label{lem:0203}  Let $\Omega$ be a bounded domain in
$\mathbb{R}^2$ and
\begin{equation}\label{0224}\Theta=\{\theta~|~\left\|\mathrm{ln}
\left(1+{\theta}\right)\right\|_{W^{1,2}(\Omega)}\leq
C_1\}.\end{equation} Then, for any $q\geq 1$, there is a
constant $C_2$ depending solely on $q$, $C_1$ and $\Omega$, such that
\begin{equation*}\label{0225}
\|\theta\|_{L^q(\Omega)}\leq C_2(q, C_1, \Omega)\quad\; \mbox{ for all }\theta\in\Theta .
\end{equation*}
\end{lem}
\begin{pf}
We use (\ref{0224}) and (\ref{0209}) to infer that
\begin{equation*}\label{0226} \Lambda :=
\left\|\mathrm{ln}\left(1+{\theta}\right)\right\|_{L_{H}(\Omega)}\leq
C_3(\Omega)\left\|\mathrm{ln}\left(1+{\theta}\right)\right\|_{W^{1,2}(\Omega)}\leq
C_3C_1.\end{equation*} Easily, it suffices to consider the case
$\Lambda\neq 0$. By the definition of the Luxemburg norm
$\|\cdot\|_{L_N(\Omega)}$, we obtain
\begin{equation*}\label{0227}
\int_\Omega
\left\{\exp{\left[\Lambda^{-1}\mathrm{ln}\left(1+{\theta}\right)\right]^2} -1\right\}
\mathrm{d}\mathbf{x}\leq 1,
\end{equation*}
which yields
\begin{equation*}\label{0228}
\int_\Omega (1+{\theta})^{\Lambda^{-2}\mathrm{ln}(1+{\theta})}
\mathrm{d}\mathbf{x}\leq 1 + |\Omega|.
\end{equation*}
Hence,
\begin{equation*}\label{0229}\begin{aligned}
\int_\Omega{\theta}^{q}\mathrm{d}x\leq &\int_{\{\theta\leq
e^{q^2}-1\}}{\theta}^{q}\mathrm{d}\mathbf{x}+\int_{\{\theta > e^{q\Lambda^2}-1\}}
{(1+ \theta)}^{q}\mathrm{d}\mathbf{x} \\
\leq &(e^{q\Lambda^2}-1)^q|\Omega|+\int_{\{\theta > e^{q\Lambda^2} -1\}}
{(1+\theta)}^{\Lambda^{-2}\mathrm{ln}(1+\theta)}\mathrm{d}\mathbf{x}\\
\leq &\left[e^{q(C_1C_3)^2}-1\right]^q|\Omega|+1+|\Omega|:=C_2(q,C_1,\Omega).
\end{aligned}\end{equation*}
This completes the proof of Lemma \ref{lem:0203}.
\hfill $\Box$
\end{pf}

Thus, with the help of Lemma \ref{0202} and the estimates
(\ref{0211})--(\ref{0215}), we conclude
\begin{equation*}\label{0230}
\|\mathrm{ln}(1+\tilde{\theta^\varepsilon}(t,|\mathbf{x}|)\|_{L^2((0,T),W^{1,2}(\Omega))}\leq
C_3(T,M_0, \Omega) ,\end{equation*} which together with Lemma
\ref{lem:0203} gives
\begin{equation}\label{0231}
\|{\theta^\varepsilon}(t,|\mathbf{x}|)\|_{L^2((0,T),L^q(\Omega^\varepsilon))}\leq
\|\tilde{\theta^\varepsilon}(t,|\mathbf{x}|)\|_{L^2((0,T),L^q(\Omega))}\leq
C_4(T, M_0, q,\Omega)\end{equation} for any $q\in [1,\infty)$, where
$\Omega^\varepsilon=\{\mathbf{x}\in {\mathbb R}^2~|~\varepsilon<|\mathbf{x}|<R\}$.
Furthermore, using (\ref{0203}), (\ref{0231}) and H\"{o}lder's inequality, we get
\begin{equation}\label{0232}
\begin{aligned}
\int_0^T\int_{\Omega^\varepsilon}|\nabla\theta^\varepsilon|
\mathrm{d}\mathbf{x}\mathrm{d}t \leq &\left(\int_0^T\int_{\Omega^\varepsilon}
\frac{|\nabla\theta^\varepsilon|^2}{(\theta^\varepsilon)^2(t,|\mathbf{x}|)}
\mathrm{d}\mathbf{x}\mathrm{d}t\right)^{1/2}\left(\int_0^T
\int_{\Omega^\varepsilon}(\theta^\varepsilon)^2 (t,|\mathbf{x}|)
\mathrm{d}\mathbf{x}\mathrm{d}t\right)^{1/2} \\[1mm]
\leq &
\left(\int_0^T\int_{\varepsilon}^R\frac{(\theta^\varepsilon_r)^2}{(\theta^\varepsilon)^2(t,r)}
r\mathrm{d}r\mathrm{d}t\right)^{1/2} C_4^{1/2}(T, M_0, 2,\Omega)\\[2mm]
\leq & C_5(T, M_0, \Omega).\end{aligned}
\end{equation}

Now, we denote
$\mathbf{u}^\varepsilon(t,\mathbf{x}):=u^\varepsilon (t,|\mathbf{x}|)\mathbf{x}/|\mathbf{x}|$,
which is a smooth function in $\Omega^\varepsilon$. By a simple calculation, we find that
\begin{equation*}  \label{0233}
|\nabla\mathbf{u}^\varepsilon|^2=(u_r^\varepsilon)^2
+\frac{(u^\varepsilon)^2}{r^2},\quad
(u_r^\varepsilon)^\zeta+\frac{u^\zeta }{r^\zeta}\leq
2^{\zeta/2}|\nabla \mathbf{u}^\varepsilon|^{\zeta} \quad\mbox{ for
any }\zeta\geq 1.\end{equation*} If we utilize (\ref{0203}),
(\ref{0231}) with $q=2$ and H\"{o}lder's inequality, we obtain
\begin{eqnarray}   \label{0234}
\int_0^T\int_{\Omega^\varepsilon}|\nabla
\mathbf{u}^\varepsilon|^{4/3}\mathrm{d}\mathbf{x}\mathrm{d}t
& \leq & \left(\int_0^T\int_{\Omega^\varepsilon} \frac{|\nabla
\mathbf{u}^\varepsilon|^2}{\theta^\varepsilon(t,|\mathbf{x}|)}
\mathrm{d}\mathbf{x}\mathrm{d}t\right)^{2/3}\left(\int_0^T
\int_{\Omega^\varepsilon} (\theta^\varepsilon )^2 (t,|\mathbf{x}|)
\mathrm{d}\mathbf{x}\mathrm{d}t\right)^{1/3}   \nonumber \\[1mm]
& \leq & \left(\frac{1}{2\mu}\right)^{2/3}
\left(\int_0^T\int_{\varepsilon}^R\frac{\mathcal{Q^\varepsilon}}{\theta^\varepsilon(t,r)}
r\mathrm{d}r\mathrm{d}t\right)^{2/3}
\left[\int_0^T\int_{\varepsilon}^R (\theta^\varepsilon)^2 (t,r)r\mathrm{d}r
\mathrm{d}t\right]^{1/3}  \nonumber \\[2mm]
& \leq & C_6(T, M_0, \Omega).
\end{eqnarray}

Since $\mathbf{u}^\varepsilon(\cdot,\mathbf{x})|_{\partial\Omega^\varepsilon}=0$,
we extend $\mathbf{u}^\varepsilon$ by zero outside
$\Omega^\varepsilon$, and employ Sobolev's inequality and (\ref{0234}) to
deduce that
\begin{equation}\label{0235}
\int_0^T\left(\int_{\Omega^\varepsilon}| \mathbf{u}^\varepsilon|^{4}
\mathrm{d}\mathbf{x}\right)^{1/3}\mathrm{d}t\leq C_{7}(T, M_0, \Omega).
\end{equation}

Finally, combining (\ref{0231})-(\ref{0234}) with (\ref{0235}), we
conclude
\begin{lem}[Global estimates of $u$ and $\theta$] \label{lem:0204}
Under the assumption of Lemma \ref{lem:0201}, there are constants $C_1$ and $C_2$, such that
\begin{eqnarray}  &&
\int_0^T\left(\int_\varepsilon^R (\theta^\varepsilon)^q (t,r)r\mathrm{d}r\right)^{2/q}
\mathrm{d}t\leq C_1(T,M_0,q,\Omega), \\
&& \int_0^T\left(\int_\varepsilon^R|\theta_r^\varepsilon(t,r)|r\mathrm{d}r\right)
\mathrm{d}t\leq C_2(T, M_0, \Omega),  \\
&&  \int_0^T\int_\varepsilon^R \left(|u^{\varepsilon
}_r|^{4/3}+\left|\frac{u^\varepsilon}{r}\right|^{4/3}\right)r\mathrm{d}r
\mathrm{d}t\leq C_2(T, M_0, \Omega),  \label{0238}  \\
&& \int_0^T\left(\int_\varepsilon^R
(u^\varepsilon)^4r\mathrm{d}r\right)^{1/3}\mathrm{d}t\leq C_2(T, M_0, \Omega).
\label{0239}
\end{eqnarray}  \end{lem}

\section{Local Estimates}
In order to taking to the limit as $\varepsilon\to 0$, we will need further uniform bounds
of higher order derivatives. Such bounds will be obtained away from the origin of
Lagrangian space in the following sense. Define a curve
$r_{h}^\varepsilon(t)$ for $h\geq 0$ by
\begin{equation}\label{0301}
h=\int_\varepsilon^{r_{h}^\varepsilon(t)}\varrho^\varepsilon(t,r)r\mathrm{d}r.
\end{equation}
Easily, by (\ref{0107}),
\begin{equation*}\frac{\partial r_{h}^\varepsilon}{\partial
t}=u^\varepsilon(t,r_{h}^\varepsilon).
\end{equation*}
Thus $r_{h}^\varepsilon(t)$ is the position at time $t$ of a fixed
fluid particle. Furthermore, an easy estimate, based on Jensen's
inequality and boundedness of
$\int_\varepsilon^R\varrho^\varepsilon\Psi(\varrho_\varepsilon^{-1})r\mathrm{d}r$
in (\ref{0203}) (see (\ref{0114})), shows
that $h\rightarrow 0$ at a uniform rate as
$r_{h}^\varepsilon(t)\rightarrow 0$. That is, given $h>0$, there is a
positive constant $C=C(h)$ independently of $\varepsilon$ and $T$, such that
\begin{equation}\label{0303}
r_{h}^\varepsilon(t)\geq C(h)^{-1} .
\end{equation}

Using (\ref{0303}), we can derive pointwise bounds for the
approximate density and temperature, which are valid away from the
origin $h=0$ of Lagrangian space, but independent of $\varepsilon$.
 The idea of deriving the pointwise boundedness was used first by
Kazhikov and Shelukhin \cite{KASVUJ}, and later adapted by Frid and Shelukhin
\cite{FHSVVS}, and by Hoff and Jenssen \cite{HDJHK} in a nontrivial
way to show a pointwise boundedness similar to that given in Lemma \ref{lem:0301} below.
 \begin{lem}[Pointwise bounds]\label{lem:0301}
Given $h>0$ and $T>0$, there is a constant $C=(T,h)$,
independent of $\varepsilon$, such that, if $r^{\varepsilon}_{h}(t)$
is given by (\ref{0301}), then
\begin{equation*}\label{0304}
C^{-1}\leq \varrho^\varepsilon(t,r)\leq C\qquad\mbox{for }r\in
[r^{\varepsilon}_{h}(t),R]\;\mbox{ and }t\in [0,T],
\end{equation*}
and
\begin{equation*}\label{0305}
\int_0^t\|\theta^\varepsilon(\tau,\cdot)\|_{h,\infty}\mathrm{d}\tau\leq
C
\end{equation*} where $\|\cdot\|_{h,\infty}$ denotes the $L^\infty$-norm
over $[r^{\varepsilon}_{h}(t),R]$.
\end{lem}
\begin{pf} Taking $n=2$ and $m=1$ in the proof of \cite[Lemma 2]{HDJHK}, we
 immediately obtain Lemma \ref{lem:0301}. \hfill $\Box$
\end{pf}

Next, we shall make use of a cut-off function which is convected with
the flow and vanishes near the origin. The cut-off function is constructed
as follows: For given $\varepsilon$ and $h$, we can fix a smooth, increasing
function $\phi_0(r)$ with $\phi_0(r)\equiv 0$ on $[0,r_h(0)]$,
$0<\phi_0(r)\leq 1$ on $(r_h(0), 2r_h(0))$ and $\phi_0(r)\equiv 1$
on $[2r_h(0),R]$, and then define $\phi(t,r)$ to be the solution of
the problem
\begin{equation}\label{0337}
\phi_t+u\phi_r=0,\qquad \phi(0,r)=\phi_0(r).
\end{equation}
We choose $\phi_0$ so that

\begin{equation*}\label{0338}
\phi_0'(r)\leq C(h)\phi_0^{(p-1)/p} \qquad\mbox{for some }p>2.
\end{equation*}

Thus, we can easily show that this boundedness persists for all time, i.e.,
\begin{equation}\label{0341}
|\phi_r(t,r(t))|\leq C(T,h)\phi(t,r(t))^{(p-1)/p}.
\end{equation}
We shall take $p$ so large that the exponent on the right-hand side of (\ref{0341})
is close to one. Notice that here we have suppressed the dependence of
$\phi$ on $\varepsilon$ and $h$.

As in \cite{HDJHK}, we now introduce three functionals of higher-order derivatives for
$(u_r^\varepsilon ,\theta^\varepsilon_r)$:
\begin{eqnarray*}  && \label{0342}
\mathcal{A}(T):=\sup_{0\leq t\leq
T}\sigma(t)\int_{r_h(t)}^R\phi(t,r)\left(u_r^\varepsilon
+\frac{u^\varepsilon}{r}\right)^2(t,r)r\mathrm{d}r  \nonumber \\[1mm]
&& \qquad\qquad
+\int_0^T\int_{r_h(t)}^R\sigma(t)\phi(t,r) (\dot{u}^\varepsilon)^2(t,r)
r\mathrm{d}r\mathrm{d}t,  \\[1mm]
&& \label{0343}
\mathcal{B}(T):=\sup_{0\leq t\leq
T}\sigma(t)\int_{r_h(t)}^R\phi(t,r) (\theta^\varepsilon)^2 (t,r) r\mathrm{d}r
+\int_0^T\int_{r_h(t)}^R\sigma(t)\phi(t,r){(\theta^\varepsilon_r)}^2(t,r)
r\mathrm{d}r\mathrm{d}t,  \\[1mm]
&&  \label{0344}
\mathcal{D}(T):=\sup_{0\leq t\leq T}
\sigma^2(t)\int_{r_h(t)}^R\phi^2(t,r)(\theta_r^\varepsilon)^2(t,r)r\mathrm{d}r
+\int_0^T\int_{r_h(t)}^R\sigma^2(t)\phi^2(t,r)(\dot{\theta}^\varepsilon)^2(t,r)
r\mathrm{d}r\mathrm{d}t, \quad\qquad
\end{eqnarray*}
where $\sigma(t)=\min\{1,t\}$, ``dot'' denotes the convective derivative
$\partial_t+u\partial_r$, and we have
again suppressed the dependence on $\varepsilon$ and $h$ for simplicity.
Thus, we have the following estimates.

\begin{lem}[Higher order boundedness]\label{lem:0303}
Let $h>0$ and  $T>0$ be given. Then there is a constant $C=C(T,h)$,
such that
\begin{equation*}\label{0345}
\int_0^T\int_{r_{h(t)}}^R\phi\left(u_r^\varepsilon+\frac{u^\varepsilon}{r}
\right)^2(t,r)r\mathrm{d}r\mathrm{d}t\leq C(T,h).
\end{equation*}
and \begin{equation*}\label{0380}\mathcal{A}(T),\ \mathcal{B}(T),\
\mathcal{D}(T)\leq C(T,h).\end{equation*}
\end{lem}
\begin{pf} This lemma can be shown following the same procedure as
in the proof of \cite[Lemmas 4 and 6]{HDJHK} with taking $m=1$, $v=0$
and $w=0$. We should point out here that in the proof one
should make use of Lemma 2.1, Lemma 3.1, (\ref{0337}) and
(\ref{0341}). \hfill$\Box$
\end{pf}

As the end of this section, we give some uniform integrability estimates. To
describe these, we define the strictly increasing, convex function $G$ by
\begin{equation*}\label{0385}\begin{aligned}
G:[1,\infty)\rightarrow [0,\infty),\quad
G(y):=y\mathrm{log}y.\end{aligned}\end{equation*} Then
$G^{-1}:[0,\infty)\rightarrow [1,\infty)$, and one can define for
$r,c>0$ the function
\begin{equation}  \label{0386}\begin{aligned}
\omega(r,c):=r+rG^{-1}\left(\frac{c}{r}\right).\end{aligned}\end{equation}
It is easy to see that for each fixed $c$ the function
$r\mapsto\omega(r,c)$ is continuous and increasing on $(0,\infty)$,
and that
\begin{equation*}\label{0387}\begin{aligned}
\lim_{r\rightarrow 0}\omega(r,c)=0.\end{aligned}\end{equation*}
Finally, if $E\subset[0,R]$, we define $|E|:=\int_Er\mathrm{d}r$.
\begin{lem}[Uniform integrability]\label{lem:0306}
Let $\omega$ be the same as in (\ref{0386}).
\begin{enumerate}[\quad \ (a)]
\item Given $b>0$ and $T>0$, there is a constant $C=C(T,b)$, such that
\begin{equation}\label{0381}
\int_{t_1}^{t_2}\left(\left\|\frac{u}{\theta^{1/2}}\right\|_{b,\infty}
+\|\mathrm{log}(\theta\vee 1)\|_{b,\infty}\right)\mathrm{d}t
\leq C(T,b)\quad\;\;\mbox{for any }t_1, t_2\in [0,T],
\end{equation}
where $\theta\vee 1=\max\{\theta,1\}$.
  \item If $\varepsilon\geq 0$, and $\varrho:[\varepsilon,R]\rightarrow
  \mathbb{R}$ is strictly positive and satisfies
\begin{equation*}\label{0388}
\left|\int_\varepsilon^R\varrho~\mathrm{log}\varrho
~r\mathrm{d}r\right|\leq C.  \end{equation*} Then, for any measure
set $E\subset [\varepsilon,R]$,
\begin{equation*}\label{0389}\begin{aligned}
\left|\int_E\varrho~r\mathrm{d}r\right|\leq
\omega(|E|,C).\end{aligned}\end{equation*}
  \item Let $b>0$ and $T>0$. Then there is a constant $C=C(T,b)$
  such that, if for $t\in [0,T]$, $E(t)$ is a measurable subset of $[b,R]$,
  and if $(\varrho,u,\theta)=(\varrho^\varepsilon,u^\varepsilon,\theta^\varepsilon)$
  is the approximate solution described at the beginning of Section 2, then
\begin{equation*}\label{0390}\begin{aligned}
\int_0^T\int_{E(t)}\varrho \theta r\mathrm{d}r\mathrm{d}t\leq
\omega\left(\int_0^T\int_{E(t)}\varrho r\mathrm{d}r\mathrm{d}t,
C(T,b)\right),\end{aligned}\end{equation*}
\begin{equation*}\label{0391}\begin{aligned}
\int_0^T\int_{E(t)}\varrho u^2 r\mathrm{d}r\mathrm{d}t\leq C(T,b)\,
\omega\left(\int_0^T\int_{E(t)}\varrho r\mathrm{d}r\mathrm{d}t,
C(T,b)\right)^{1/4}.  \end{aligned}\end{equation*}
\end{enumerate}
\end{lem}
\begin{pf} The estimate (\ref{0381}) is a consequence of the entropy estimate (\ref{0203}).
The uniform integrability bounds in (b) and (c), which are important
in showing the limit of these approximate solutions to be a weak
solution in Section 4, can be obtained using arguments similar to
those used in \cite[Lemmas 8 and 9]{HDJHK}, and hence we omit the
proof here.     \hfill$\Box$
\end{pf}

\section{Proof of Theorem \ref{thm:0101}}
By virtue of the a priori estimates derived in Sections 2 and 3, we are now able to
prove our main theorem by taking appropriate limits in a manner
analogous to that in \cite{HDJHK}.

To begin with, we denote by $H_\delta$ a standard mollifier (in $r$)
of width $\delta$, and for $\varepsilon>\delta$ we define the smooth
approximate initial data
$(\varrho_0^{\varepsilon,\delta},{u}_0^{\varepsilon,\delta},\theta_0^{\varepsilon,\delta}$)
to $(\varrho_0,{u}_0,\theta_0)$ as follows:

1) Extend $\varrho_0$ by its average value outside
  $[\varepsilon,R]$, mollify with $H_\delta$, restrict to
  $[\varepsilon,R]$, and then multiply by a constant to normalize
  the total mass to be $M_0=\int_0^R\varrho_0r\mathrm{d}r$. The
  resulting density function is denoted by $\varrho_0^{\varepsilon,\delta}(r)$.

  2) Redefine $u_0$ to be zero on $[0,2\varepsilon]$ and $[R-2\delta,
  R]$, then mollify with $H_\delta$ to get $u_0^{\varepsilon,\delta}$.
  Note that $u_0^{\varepsilon,\delta}$ is identically zero in a neighborhood of
  $r=\varepsilon$ or $R$.

 3) Redefine $\theta_0$ to be its average value on $[0,2\varepsilon]$ and
  $[R-2\delta]$, then mollify with $H_\delta$ to get $\theta_0^{\varepsilon,\delta}$.
  Note that $\theta_0^{\varepsilon,\delta}$ is constant in a neighborhood
  of $r=\varepsilon$ or $R$.

The resulting data $\varrho_0^{\varepsilon,\delta}$,
$u_0^{\varepsilon,\delta}$, $\theta_0^{\varepsilon,\delta}$ then
satisfy the hypotheses (\ref{0112}) and (\ref{0113}) with the
constants independent of $\varepsilon$ and $\delta$. Thus, there is
a global-in-time smooth solution
($\varrho^{\varepsilon,\delta},{u}^{\varepsilon,\delta},\theta^{\varepsilon,\delta}$)
of the system (\ref{0107})--(\ref{0109}) with the initial boundary
conditions (\ref{0109b}) and (\ref{j1}). This is a result of Frid
and Shelukhin's work \cite{FHSVVS} in the annular domain
$(\varepsilon ,R)$. Next, we want to pass to the limit to get a
global weak solution. As in (\ref{0301}), we define the particle
path $r_h^{\varepsilon,\delta}(t)$ associated with the approximate
solution
($\varrho^{\varepsilon,\delta},{u}^{\varepsilon,\delta},\theta^{\varepsilon,\delta}$)
by \begin{equation}\label{0401}
h=\int_\varepsilon^{r_h^{\varepsilon,\delta}(t)}\varrho^{\varepsilon,\delta}(t,r)r\mathrm{d}r,
\qquad h, \varepsilon, \delta >0.
\end{equation}

\subsection{Convergence of the Approximate Solutions}
By the a priori estimates established in Lemmas
\ref{lem:0301}--\ref{lem:0303}, \ref{lem:0201} and \ref{lem:0204},
we have the following three propositions, which imply that there is
a subsequence $(\varepsilon_j,\delta_j)\rightarrow (0,0)$, such that
the approximate solutions and their associated particle paths are
convergent.
\begin{pro}\label{pro:0401}
Let
$(\varrho^{\varepsilon,\delta},{u}^{\varepsilon,\delta},\theta^{\varepsilon,\delta}$)
and $r_h^{\varepsilon,\delta}(t)$ be as described above.
\begin{enumerate}[\quad \ (a)]
  \item There is a subsequence $(\varepsilon_j,\delta_j)\to (0,0)$,
   such that $r_h^{\varepsilon,\delta}(t)$ converges uniformly
  for ($t, h$) in any compact subset of $[0,\infty)\times (0,M_0)$, and
  the limit $r_h(t)$ is H\"{o}lder-continuous in ($t,h$) on any compact set.
  \item If $\underline{r}(t):=\lim_{h\to 0} r_h(t)$, then
$\underline{r}(t)\in C^{0,1/4}_{\mathrm{loc}}[0,\infty)$ and $\lim_{t\to 0}\underline{r}(t)=0$.
  \item If the ``fluid region" $\mathcal{F}$ is defined by
\begin{equation*}\label{0402}
\begin{aligned}
\mathcal{F}:=\{(t,r)~|~\underline{r}(t)<r\leq R,\ 0\leq t<\infty\},
\end{aligned}
\end{equation*}
then $\mathcal{F}\cap \{t>0\}\cap \{r<R\}$ is an open set.
\end{enumerate}
\end{pro}
\begin{pro}\label{pro:0402}
Let the hypotheses of Proposition \ref{pro:0401} be satisfied.
Then there is a further subsequence, still denoted by
$(\varepsilon_j,\delta_j)$, and limiting functions $u$ and $\theta$, such that
\begin{equation*}\label{0405}
\begin{aligned}
u^{\varepsilon_j, \delta_j}\rightarrow u,\ \theta^{\varepsilon_j,
\delta_j}\rightarrow \theta
\end{aligned}
\end{equation*}
uniformly on any compact subset of $\mathcal{F}\cap \{t>0\}$. The
functions $u$ and $\theta$ are H\"{o}lder-continuous on any
compact set. Furthermore, for any $T>0$,
\begin{equation}\label{0406}
\begin{aligned}
u^{\varepsilon_j, \delta_j}\rightarrow u\mbox{ weakly in }L^{4/3}
((0,T),W^{1,4/3}_{\mathrm{loc}}(0,R])
\end{aligned}
\end{equation}
\end{pro}
\begin{pro}
Assume that the hypotheses of Proposition 4.2 hold. Then there is a
further subsequence ($\varepsilon_j$, $\delta_j)\rightarrow$ ($0$,
$0$) and a function $\varrho(t,r)$ such that
\begin{equation*}\label{pro:0403}
\begin{aligned}
\varrho^{\varepsilon_j, \delta_j}(t,\cdot)\rightarrow \varrho
(t,\cdot)\mbox{ in }H^{-1}([\underline{r}(t)+\eta,R],r\mathrm{d}r)
\end{aligned}
\end{equation*}
and\begin{equation*}\label{0411}
\begin{aligned}
\varrho^{\varepsilon_j, \delta_j}(t,\cdot)\rightharpoonup \varrho
(t,\cdot)\mbox{ in }L^{2}([\underline{r}(t)+\eta,R],r\mathrm{d}r)
\end{aligned}
\end{equation*}
for all $t\in [0,T]$ and all $\eta>0$. In addition, if
$\varrho^{\varepsilon_j, \delta_j}(t,\cdot)$ is taken to be zero for
$r\leq \varepsilon_j$, then
\begin{equation*} \label{0412}
\varrho^{\varepsilon_j, \delta_j}(t,\cdot)\rightarrow 0\mbox{ in }L^{1}
([0,\underline{r}(t)],r\mathrm{d}r)\qquad \mbox{when }\underline{r}(t)>0.
\end{equation*}
Also, for $h>0$ and $T>0$, there is constant $C=C(T,h)$, such that
\begin{equation*}\label{0413}
C^{-1}(T,h)\leq \varrho \leq C(T,h)\quad\;\;\mbox{for }0\leq t\leq T
\mbox{ and }r_h(t)\leq r\leq R.
\end{equation*}
Finally, for $h>0$ and $t\geq 0$,
\begin{equation}\label{0414}
\begin{aligned}
h=\int_{\underline{r}(t)}^{r_h(t)}\varrho(t,r) r\mathrm{d}r.
\end{aligned}
\end{equation}
\end{pro}

By virtue of the a priori estimates established in Sections 2 and 3, one can
show Propositions 4.1--4.3 in the same manner as that in the proof of
 \cite[Propositions 1--3]{HDJHK}, except (\ref{0406}) of Proposition 4.2
which are obtained by applying the uniform global estimates given in
Lemma \ref{lem:0204}. In addition, the identity (\ref{0414}) shows
that mass is conserved for the limiting solution
\begin{equation}\label{0428}
\begin{aligned}
M_0=\int_{\underline{r}(t)}^R\varrho(t,r)
r\mathrm{d}r=\int_0^R\varrho_0(r)r\mathrm{d}r.\end{aligned}\end{equation}

\subsection{Weak Forms of the Navier-Stokes Equations}
We now turn to the proof that the limiting functions are
indeed a weak solution of the Navier-Stokes equations in
$[0,\infty)\times \Omega$ in the sense of Theorem 1.1.

First, the limiting functions $\varrho$, $u$ and $\theta$ have been
defined in the fluid region $\mathcal{F}$ but not elsewhere. We
therefore define $\varrho$, $\varrho u$ and $\varrho\theta$ to be
identically zero in the vacuum region $\mathcal{F}^c$. As in Section
1 we let $r(\mathbf{x})=|\mathbf{x}|$ and define the velocity vector
$\mathbf{u}:[0,\infty)\times \bar{\Omega}$ by
\begin{equation}\label{0429}
\mathbf{u}(t,\mathbf{x})=u(t,r)\frac{\mathbf{x}}{r}.
\end{equation}

Abusing notation slightly, we also write $\varrho(t,\mathbf{x})$ and
$\theta(t,\mathbf{x})$ in place of $\varrho(t,r(\mathbf{x}))$ and
$\theta(t,r(\mathbf{x}))$. Similar notation applies to the
approximate solutions, for which we now write $u^j$ in place of
$u^{\varepsilon_j,\delta_j}$, etc.

We first show that ($\varrho$, $\mathbf{u}$, $\theta$) satisfies
the weak form (\ref{0116}) of the mass equation.

\begin{pro}\label{pro:0404}
Let ($\varrho$, $\mathbf{u}$, $\theta$) be the limit described above
in Propositions \ref{pro:0401}--4.3. Then,
\begin{enumerate}[\quad \ (a)]
  \item The weak form (\ref{0116}) of the mass equation holds for any $C^1$
  test function $\phi:[t_1,t_2]\times \bar{\Omega}\rightarrow
  \mathbb{R}$.
  \item $\varrho\in C([0,\infty),W^{1,\infty}(\Omega))^*$.
  \item $\varrho^{1/2}\mathbf{u}\in L^\infty([0,\infty),L^2(\Omega))$.
 \item $\mathbf{u}\in L^{4/3}(0,T;W^{1,4/3}(\Omega))$.
\end{enumerate}\end{pro}
\begin{pf} The assertions (a)--(c) follow from the analogous arguments as in the proof
of \cite[Proposition 4]{HDJHK}, where we have made use of Lemma 3.3
and Proposition \ref{pro:0401}-\ref{pro:0403}.

To show (d), we use (\ref{0238}), (\ref{0239}), (\ref{0406}) and the
lower semi-continuity to deduce that
\begin{equation}\label{0443}
\int_0^T\int_0^R
\left(|u_r|^{4/3}+\left|\frac{u}{r}\right|^{4/3}\right)r\mathrm{d}r\mathrm{d}t
\leq \lim_{b\to 0}\lim_{j\to\infty}\int_0^T\int_b^R \left(|u^{j}_r|^{4/3}
+\left|\frac{u^j}{r}\right|^{4/3}\right)r\mathrm{d}r\mathrm{d}t \leq C(T)
\end{equation} and
\begin{equation*} \label{0444}
\int_0^T\left(\int_0^R u^4r\mathrm{d}r\right)^{1/3}\mathrm{d}t\leq
\lim_{b\to 0}\lim_{j\to\infty}\int_0^T\left (\int_b^R u_j^4r\mathrm{d}r\right)^{1/3}
\mathrm{d}t\leq C(T).
\end{equation*}

We can compute that
\begin{equation} \label{04444}\partial_{x_j}\left[u(t,|\mathbf{x}|)\frac{x_i}{|\mathbf{x}|}\right]
=u_r(t,|\mathbf{x}|)\frac{x_ix_j}{|\mathbf{x}|^2}+u(t,|\mathbf{x}|)\left(\frac{\delta_{ij}}{|\mathbf{x}|}
-\frac{x_ix_j}{|\mathbf{x}|^3}\right),\ 1\leq i,\ j\leq 2,
\end{equation}
thus, $|\nabla \mathbf{u}|^{4/3}\leq
2^{2/3}\left(|u_r|^{4/3}+|u/r|^{4/3}\right)$ and $|
\mathbf{u}|^2=u^2$. Hence, we find that $\mathbf{u}\in
L^{4/3}(0,T;$$W^{1,4/3}$ $(\Omega))$. This completes the proof.
\hfill $\Box$
\end{pf}

Next, we show that the weak form of the momentum equations in the spherically symmetric case holds.
\begin{lem}\label{lem:0401}
Let $\varrho$, $u$ and $\theta$ be the functions given in
Propositions 4.2 and 4.3. Let $t_1<t_2$ and $\phi$ be a
$C^1$-function on $[t_1,t_2]\times [0,R]$, such that
$\phi(t,0)=\phi(t,R)=0$ for $t\in [t_1,t_2]$. Then, the following
identity holds.
\begin{equation}\label{0445}\begin{aligned}
&\int_0^R\varrho u\phi
r\mathrm{d}r\bigg|_{t_1}^{t_2}-\int_{t_1}^{t_2}\int_0^R\left[\varrho
u\phi_t+{\varrho
u^2}\phi_r+P(\varrho,\theta)\left(\phi_r+\frac{\phi}{r}\right)\right]r\mathrm{d}r\mathrm{d}t\\
=&\int_{t_1}^{t_2}\int_0^R \varrho f\phi
r\mathrm{d}r\mathrm{d}t-\nu\int_{t_1}^{t_2}\int_0^R\left(u_r+\frac{u}{r}\right)
\left(\phi_r+\frac{\phi}{r}\right)r\mathrm{d}r\mathrm{d}t.\end{aligned}
\end{equation}
\end{lem}
\begin{pf}
We first consider a simpler case in which the test function vanishes
in a neighborhood of the origin. Assume that $\psi$ is a
$C^1$-function on $[t_1,t_2]\times [0,R]$ satisfying $\psi\equiv 0$
on $[0,b]$ for some $b>0$. Then, applying the proof of \cite[Lemma
10 (a)]{HDJHK} and combining with the weak convergence (\ref{0406}),
we can easily show that the weak form of the momentum equation
(\ref{0108}) holds for the test function $\psi$:
\begin{equation}\label{0446} \begin{aligned}
& \int_0^R\varrho u\psi
r\mathrm{d}r\bigg|_{t_1}^{t_2}-\int_{t_1}^{t_2}\int_0^R\left[\varrho
u\psi_t+{\varrho
u^2}\psi_r+P\left(\psi_r+\frac{\psi}{r}\right)\right]r\mathrm{d}r\mathrm{d}t\\
=& \int_{t_1}^{t_2}\int_0^R \varrho f\psi
r\mathrm{d}r\mathrm{d}t-\nu\int_{t_1}^{t_2}\int_0^R\left(u_r+\frac{u}{r}\right)
\left(\psi_r+\frac{\psi}{r}\right)r\mathrm{d}r\mathrm{d}t.\end{aligned}
\end{equation}

To extend the identity (\ref{0446}) to the case that test functions
are supported in $[0,R]$, we fix an increasing $C^1$ function
$\chi:[0,\infty)\rightarrow [0,1]$ with $\chi\equiv 0$ on $[0,1]$
and $\chi\equiv 1$ on $[2,\infty)$, and define
$\chi^b(r):=\chi(r/b)$ for $b>0$. Let $\phi$ be a $C^1$ function on
$[t_1,t_2]\times [0,R]$ such that $\phi(t,0)=\phi(t,R)=0$ for $t\in
[t_1,t_2]$, and define $\phi^b:=\chi^b\phi$. Then, the previous
lemma applies to the test functions $\phi^b=\chi^b\phi$. We obtain
\begin{equation}\label{0451}     \begin{aligned}
& \int_0^R\varrho u\chi^b\phi
r\mathrm{d}r\bigg|_{t_1}^{t_2}-\int_{t_1}^{t_2}\int_0^R\bigg[\varrho
u\chi^b\phi_t+{\varrho u^2}(\chi^b\phi)_r+P(\chi^b\phi)_r
+\frac{P\chi^b\phi}{r}\bigg]r\mathrm{d}r\mathrm{d}t\\
= & \int_{t_1}^{t_2}\int_0^R\varrho f\chi^b\phi r
\mathrm{d}r\mathrm{d}t-\nu\int_{t_1}^{t_2}\int_0^R
\left(u_r+\frac{u}{r}\right)\left[(\chi^b\phi)_r+\frac{\chi^b\phi}{r}
\right]_rr\mathrm{d}r\mathrm{d}t.   \end{aligned}
\end{equation}
The first, second, fifth and sixth terms in (\ref{0451}) converge to the
corresponding terms in (\ref{0445}) as $b\to 0$ by virtue of the
dominated convergence theorem. For the third term we have
\begin{equation}\label{0452}
\int_{t_1}^{t_2}\int_0^R{\varrho u^2}(\chi^b\phi)_rr\mathrm{d}r\mathrm{d}t=
\int_{t_1}^{t_2}\int_0^R{\varrho u^2}(\chi^b_r\phi+
\chi^b\phi_r)r\mathrm{d}r\mathrm{d}t ,
\end{equation}
and the second term on the right-hand side of (\ref{0451}) clearly
tends to the third term in (\ref{0445}) as $b\to 0$. Since
$\phi(t,0)=0$, we can write $\phi(t,r)=r\varphi(t,r)$ for some
smooth $\varphi$. Thus, due to $|\partial_r\chi^b|\leq C/b$ we can
bound the first term on the right-hand side of (\ref{0452}) by
\begin{equation}\label{0453}
\int_{t_1}^{t_2}\int_0^b{\varrho u^2}
C\frac{r}{b}\left|\varphi(r)\right|r\mathrm{d}r\mathrm{d}t\leq
C\int_{t_1}^{t_2}\int_b^{2b}{\varrho u^2}r\mathrm{d}r\mathrm{d}t ,
\end{equation}
which goes to zero as $b\to 0$ by utilizing the boundedness of the limiting
energy. Moreover, the same argument applies to the fourth term in (\ref{0451}).
Finally, for the last term on the right-hand side of (\ref{0451}),
we have
\begin{equation} \label{0454}
\int_{t_1}^{t_2}\int_0^R \left(u_r+\frac{u}{r}\right)
\left(\phi_r^b+\frac{\phi^b}{r}\right)r\mathrm{d}r\mathrm{d}t\\
=\int_{t_1}^{t_2}\int_0^R
\left(u_r+\frac{u}{r}\right)\left(\phi\chi_{r}^b+\phi_r\chi^b
+\frac{\chi^b\phi}{r}\right)r\mathrm{d}r\mathrm{d}t.
\end{equation}
Similarly to (\ref{0453}), we deduce that
\begin{equation}\label{0455}
\left|\int_{t_1}^{t_2}\int_0^R \left(u_r+\frac{u}{r}\right)
\phi\chi_{r}^br\mathrm{d}r\mathrm{d}t\right| \leq C\int_{t_1}^{t_2}\int_{b}^{2b}
\left(u_r+\frac{u}{r}\right)r\mathrm{d}r\mathrm{d}t ,
\end{equation}
which tends to zero as $b\to 0$ by (\ref{0443}). On the
other hand, letting $b\to 0$ in (\ref{0454}), using the dominated
convergence theorem and (\ref{0455}), we conclude
\begin{equation*}\label{0456}
\int_{t_1}^{t_2}\int_0^R \left(u_r+\frac{u}{r}\right)
\left(\phi_{r}\chi^b
+\frac{\phi\chi^b}{r}\right)r\mathrm{d}r\mathrm{d}t
\to\int_{t_1}^{t_2}\int_0^b
\left(u_r+\frac{u}{r}\right)\left(\phi_{r}+\frac{\phi}{r}\right)r\mathrm{d}r\mathrm{d}t.
\end{equation*}
This completes the proof of the lemma.
\end{pf}
\begin{rem}
Applying (\ref{0238})-(\ref{0239}) and the proof of  Lemma
\ref{lem:0401} to \cite[Lemma 11]{HDJHK}, we can see that
\begin{equation*}
\lim_{R\rightarrow 0}\lim_{j\rightarrow
\infty}\mathcal{U}(j,\phi^R)=\nu\int_{t_1}^{t_2}\int_0^b\left(u_r+\frac{u}{r}\right)
\left(\phi_r+\frac{\phi}{r}\right)r\mathrm{d}r\mathrm{d}t
\end{equation*}
still holds for cylindrically symmetric case in \cite[Lemma
11]{HDJHK}.
\end{rem}

Now, we are able to show that the weak form (\ref{0117}) for the momentum equations
in Cartesian coordinates is satisfied.
\begin{pro}\label{pro:0405} The weak form (\ref{0117}) of the momentum equations,
as stated in Theorem 1.1 (d), holds.

\end{pro}
\begin{pf}
Given $\psi$ as described in the theorem, we define
\begin{equation*}\label{0457}
\begin{aligned}\phi(t,r):=\int_S\psi(t,r\mathbf{y})y_i\mathrm{d}S_\mathbf{y}\mbox{ for fixed }i=1,\ 2,
\end{aligned}
\end{equation*}
where $S\subset \mathbb{R}^2$ denotes the unit circle. Then,
$\phi(t,0)=\phi(t,R)\equiv 0$. It thus follows from Lemma \ref{lem:0401} that
\begin{equation}\label{0458}  \begin{aligned}
&\int_0^R\varrho u\phi
r\mathrm{d}r\bigg|_{t_1}^{t_2}-\int_{t_1}^{t_2}\int_0^R\left[\varrho
u\phi_t+{\varrho u^2}\phi_r+P\left(\phi_r+\frac{\phi}{r}\right)\right]r\mathrm{d}r\mathrm{d}t\\
=&\int_{t_1}^{t_2}\int_0^R \varrho f\phi
r\mathrm{d}r\mathrm{d}t-\nu\int_{t_1}^{t_2}\int_0^R\left(u_r+\frac{u}{r}\right)
\left(\phi_r+\frac{\phi}{r}\right)r\mathrm{d}r\mathrm{d}t.\end{aligned}
\end{equation}
We convert each of the terms in (\ref{0458}) to integrals in Cartesian coordinates
involving $\psi$. The treatment of the terms involving derivations are very much similar
to those in the proof of Proposition 4.4 (see \cite[Proposition 4]{HDJHK}),
except for the last term, which we deal with in details. We may
rewrite the last term of (\ref{0458}) as
\begin{equation}\label{0459}
\begin{aligned}& \int_{t_1}^{t_2}\int_{0}^R\left(u_r+\frac{u}{r}\right)
\left(\phi_r+\frac{\phi}{r}\right)r\mathrm{d}r\mathrm{d}t\\
& \quad
=\int_{t_1}^{t_2}\int_{0}^R\left[\left(\frac{u}{r}\right)_r\phi_rr^2
+\frac{u}{r}(r\phi)_r\right]\mathrm{d}r\mathrm{d}t  \\
& \quad =\int_{t_1}^{t_2}\int_{0}^Rr\left(\frac{u}{r}\right)_r\left(\int_S\psi_{x_k}(t,r\mathbf{y})y_i{y_k}
\mathrm{d}S_\mathbf{y}\right)r\mathrm{d}r\mathrm{d}t\\
& \qquad +\int_{t_1}^{t_2}\int_{0}^R\frac{u}{r}\left(\int_S\psi_{{x_i}}(t,r\mathbf{y})
\mathrm{d}S_\mathbf{y}\right)r\mathrm{d}r\mathrm{d}t\\
& \quad =\int_{t_1}^{t_2}\int_\Omega\left[\left(\frac{u(t,|\mathbf{x}|)}{|\mathbf{x}|}\right)_r
\frac{x_k x_i}{|\mathbf{x}|}+\frac{u}{|\mathbf{x}|}\delta_{ik}\right]
\psi_{x_k}(t,\mathbf{x})r\mathrm{d}\mathbf{x}\mathrm{d}t\\
& \quad  = \int_{t_1}^{t_2}\int_\Omega\left(u(t,|\mathbf{x}|)\frac{x_i}{|\mathbf{x}|}\right)_{x_k}
\psi_{x_k}(t,r\mathbf{y})r\mathrm{d}\mathbf{x}\mathrm{d}t\\
& \quad =\int_{t_1}^{t_2}\int_\Omega \nabla
u_i(t,|\mathbf{x}|)\nabla\psi\mathrm{d}\mathbf{x}\mathrm{d}t,
\end{aligned}
\end{equation}
where the repeated indexes should be summed. Notice that
$\partial_{x_j}u_i(t,|\mathbf{x}|)=\partial_{x_i}u_j(t,|\mathbf{x}|)$
(see (\ref{04444})), we have
\begin{equation*}\label{0460}
\begin{aligned}\int_{t_1}^{t_2}\int_\Omega \nabla u_i(t,|\mathbf{x}|)\cdot\nabla
\psi\mathrm{d}\mathbf{x}\mathrm{d}t=&\int_\Omega
\partial_{x_k}u_i(t,|\mathbf{x}|)\partial_{x_k}\psi\mathrm{d}\mathbf{x}\mathrm{d}t
\\=&\int_{t_1}^{t_2}\int_\Omega
\partial_{x_i}u_k(t,|\mathbf{x}|)\partial_{x_k}\psi\mathrm{d}\mathbf{x}\mathrm{d}t
\\=&\int_{t_1}^{t_2}\int_\Omega
\mathbf{u}_{x_i}\cdot\nabla \psi\mathrm{d}\mathbf{x}\mathrm{d}t ,
\end{aligned}
\end{equation*}
whence, by recalling $\nu=\lambda+2\mu$, the last term on the right-hand of
(\ref{0458}) can be written
\begin{equation*}\label{0461}
(\lambda+\mu)\int_{t_1}^{t_2}\int_\Omega
\mathbf{u}_{x_i}\cdot\nabla\psi \mathrm{d}\mathbf{x}
\mathrm{d}t+\mu\int_{t_1}^{t_2}\int_\Omega \nabla u_i\nabla
\psi\mathrm{d}\mathbf{x}\mathrm{d}t.
\end{equation*}
This completes the proof of Proposition 4.5. \hfill $\Box$
\end{pf}
\begin{pro}
The weak form (\ref{0118}) of the energy equation, as stated in (e)
of Theorem \ref{thm:0101}, holds. The total energy $\mathcal{E}$ of
the the limiting functions satisfies (f) of Theorem \ref{thm:0101}.
\end{pro}
\begin{pf}
The proof is the same as that of \cite[Propositions 7 and 8]{HDJHK}, and hence we omit
the proof here.\hfill $\Box$
\end{pf}

Thus, we have completed the proof of Theorem 1.1. In fact, Part (a)
of Theorem \ref{thm:0101} is just Proposition \ref{pro:0401} with
the semicontinuity used implicitly in the proof. The existence and
regularity of $\varrho$, $\mathbf{u}$, $\theta$ asserted in (b) of
Theorem \ref{thm:0101} follow from Propositions \ref{pro:0402} and
4.3. The weak forms of the mass and momentum equations are proved in
Propositions \ref{pro:0404} and \ref{pro:0405}, from which the
regularity assertions in (c) and (d) of Theorem \ref{thm:0101}
follow immediately. Finally, the results in (e) and (f) of the
Theorem \ref{thm:0101} are proved by applying Proposition 4.6.

%

\end{document}